\documentclass[11pt,reqno,oneside]{amsart}
\usepackage{verbatim,amsmath,amssymb,cite,epsfig,xspace}

\numberwithin{equation}{section}

\theoremstyle{plain}
\newtheorem{theorem}{Theorem}[section]

\theoremstyle{definition}

\theoremstyle{remark}
\newtheorem{remark}{Remark}[section]

\newcommand{\ignore}[1]{}

\newcommand{\yb}{\mathbf y}

\newcommand{\F}{{\mathcal F}}

\newcommand{\E}{{\mathcal E}}

\newcommand{\rhoa}{\bar\rho^a_\ell(\mathbf{y})}

\newcommand{\rhoc}{\bar\rho^c_\ell(\mathbf{y})}
\newcommand{\rhocc}{\bar\rho^c_{\ell+1}(\mathbf{y})}
\newcommand{\rhoqnl}{\bar\rho^{qnl}_\ell(\mathbf{y})}
\newcommand{\la}{\langle}
\newcommand{\ra}{\rangle}

\newcommand{\half}{\frac{1}{2}}

\newcommand{\eps}{\epsilon}

\begin{document}

\title
[Quasi-Nonlocal Quasicontinuum Approximation of the Embedded Atom Model]
{An Analysis of the Quasi-Nonlocal Quasicontinuum Approximation of the Embedded Atom Model}
\author{Xingjie Helen Li and Mitchell Luskin}


\thanks{
This work was supported in part by DMS-0757355,
 DMS-0811039,  the Institute for Mathematics and
Its Applications, and
 the University of Minnesota Supercomputing Institute.
  This work was also supported by the Department of Energy under Award Number
DE-SC0002085.
}



\date{\today}

\begin{abstract}
The quasi-nonlocal quasicontinuum method (QNL) is a consistent
hybrid coupling method for atomistic and continuum models.
Embedded atom models are empirical many-body potentials
that are widely used for FCC metals such as copper and aluminum.
In this paper, we consider the QNL method for EAM potentials,
and we give a stability and error
analysis for a chain with next-nearest neighbor interactions.
We identify conditions for the pair potential, electron density
function, and embedding function so that the lattice stability of the atomistic and the
EAM-QNL models are asymptotically equal.
\end{abstract}

\maketitle{\thispagestyle{empty} \maketitle
\section{Introduction}
Hybrid atomistic-to-continuum methods couple atomistic regions surrounding
defects with continuum regions to achieve the accuracy of
the atomistic model and the efficiency of the continuum model.
Quasicontinuum hybrid methods utilize the Cauchy--Born rule for the
energy density in the continuum region~\cite{Ortiz:1995a}. The original
quasicontinuum energy \cite{Ortiz:1995a}
 (denoted QCE) has
interfacial forces (called ``ghost forces'') for a uniform strain
~\cite{Shenoy:1999a,Dobson:2008a}.
Thus, uniform strain is not an equilibrium solution for
the QCE energy (even though uniform strain is an equilibrium for
 purely atomistic and for purely coarse-grained continuum models).

More accurate atomistic-to-continuum coupling methods have been
proposed to remedy the QCE model. The ghost force correction method
(GFC) achieves an increased accuracy by adding a correction to the
ghost forces as a dead load during a quasistatic
process~\cite{Shenoy:1999a,Dobson:2008a,dobsonluskin08,Miller:2008,qcf.iterative}.
The GFC method can be viewed as a stationary iterative method
~\cite{Dobson:2008a,dobsonluskin08,Miller:2008,qcf.iterative} to solve
the force-based quasicontinuum aproximation (QCF) using QCE as a preconditioner.
More accurate coupling can be achieved by using a more accurate preconditioner
or by using GMRES acceleration to solve the QCF
equilibrium equations~\cite{doblusort:qcf.stab,dobs-qcf2,curt03,Miller:2003a},
but the non-conservative and indefinite QCF equilibrium equations make the iterative
solution and the determination of lattice stability more
challenging~\cite{qcf.iterative}.

An alternative approach is to develop a quasicontinuum energy that
is more accurate than QCE. We will call a QC energy {\it consistent}
if it does not have ghost forces for a uniformly strained lattice.
The quasi-nonlocal energy (QNL) was the
first consistent quasicontinuum energy~\cite{Shimokawa:2004}.
For a one dimensional chain, the
original QNL method is restricted to next-nearest neighbor
interactions~\cite{Shimokawa:2004}. The QNL method
for pair interaction potentials was extended to
finite range interactions in ~\cite{LuskinXingjie} and to
two dimensional finite range problems in ~\cite{shapeev}.

In this paper, we formulate a one-dimensional QNL
energy for the embedded atom model (EAM) following
~\cite{Shimokawa:2004}.
The embedded atom model~\cite{Foiles,Mishin,Johnson:1989} is an empirical many-body
potential that is widely used to model FCC metals such as copper
and aluminum.
We then give an analysis of the stability and error
for the EAM-QNL approximation in the next-nearest neighbor case for a periodic chain.

We identify conditions for the pair potential, electron density
function, and embedding function so that the lattice stability of the atomistic and the
EAM-QNL models are asymptotically equal.
We also show in Remark~\ref{embed} that the
atomistic and EAM-QNL models can be less stable than the local quasicontinuum model (EAM-QCL),
which is the EAM-QNL model with no atomistic region,
if the above conditions on the pair potential, electron density
function, and embedding function are not satisfied.

Many theoretical analyses of QC models have been given based on
pair-potential interactions
\cite{Dobson:2008b,BadiaParksBochevGunzburgerLehoucq:2007,E:2006,doblusort:qce.stab,mingyang,LinP:2006a,LuskinXingjie}.
In this paper, we give an analysis of the stability and
accuracy of a linearization of the quasi-nonlocal method  for the
EAM potential in one dimension with next-nearest neighbor
interactions.  A nonlinear {\it a priori} and {\it
a posteriori} error analysis for the QNL model with next-
nearest-neighbor pair potential interaction in one dimension was
given in ~\cite{ortner:qnl1d}.  We think that a similar nonlinear analysis
using the inverse function theorem can
be done for this model, but we restrict our presentation in this paper to
the linear analysis for simplicity.

In Section~\ref{notation}, we present the notation used in this
paper. We define the displacement space $\mathcal{U}$ and the
deformation space $\mathcal{Y}_{F}$. We then introduce the norms we
will use to estimate the modeling error and the displacement
gradient error.
In Section~\ref{qnl}, we introduce the QNL model with next-nearest neighbor
interaction for the EAM potential.

In Section~\ref{sharp}, we give sharp stability estimates for both
the fully atomistic model and the EAM-QNL model for a
uniformly strained chain. Sharp stability estimates are necessary to determine
whether quasicontinuum methods (or other coupling methods) are
accurate near instabilities such as defect formation or crack
propagation~\cite{doblusort:qce.stab,doblusort:qcf.stab}. Similar
stability estimates for the fully atomistic and fully local
quasi-continuum (QCL) models can also be obtained by discrete Fourier
analysis~\cite{HudsOrt:a}.

In section~\ref{convergence}, we study the convergence rate of the
EAM-QNL model. We compare the equilibrium solution of the EAM-QNL model
with that of the fully atomistic model, and we use the negative norm
estimation method \cite{dobs-qcf2,LuskinXingjie} to obtain an
optimal rate of convergence of the strain error.  The error estimate
depends only on the smoothness of the strain in the continuum region
and holds near lattice instabilities, thus demonstrating that the
QNL method for the EAM potential can give a small error if defects
are captured in the atomistic region.


\section{Notation}
\label{notation}
In this section, we
present the notation used in this paper. We define the scaled reference
lattice
\[
\eps \mathbb{Z}:= \{\eps\ell : \ell\in\mathbb{Z}\},
\]
 where $\eps >0$
scales the reference atomic spacing and $\mathbb{Z}$ is the set of
integers.  We then deform the reference lattice $\eps \mathbb{Z}$
uniformly into the lattice
\[
F\eps \mathbb{Z}:= \{F\eps\ell : \ell\in\mathbb{Z}\}
\]
where $F >0$ is the macroscopic deformation gradient, and we define
the corresponding deformation $\mathbf{y}_{F}$ by
\[
(\mathbf{y}_{F})_{\ell}:=F\eps \ell  \quad \text{for } -\infty
<\ell<\infty.
\]
For simplicity, we consider the space $\mathcal{U}$ of $2N$-periodic
zero mean displacements $\mathbf{u}=(u_{\ell})_{\ell \in
\mathbb{Z}}$ from $\mathbf{y}_{F}$ given by
\[
\mathcal{U}:=\bigg\{\mathbf{u} : u_{\ell+2N}=u_{\ell}
\text{ for }\ell\in \mathbb{Z},\,
\text{and}\sum_{\ell=-N+1}^{N}u_{\ell}=0\bigg\},
\]
and we thus admit deformations $\mathbf{y}$ from the space
\[
\mathcal{Y}_{F}:=\{\mathbf{y}:
\mathbf{y}=\mathbf{y}_{F}+\mathbf{u}\text{ for some }\mathbf{u}\in
\mathcal{U}\}.
\]
We set $\eps=1/N$ throughout so that the reference length of the
periodic domain is fixed.

 We
define the discrete differentiation operator, $D\mathbf{u}$, on
periodic displacements by
\[
(D\mathbf{u})_{\ell}:=\frac{u_{\ell}-u_{\ell-1}}{\epsilon}, \quad
-\infty<\ell<\infty.
\]
We note that $\left(D\mathbf{u}\right)_{\ell}$ is also $2N$-periodic
in $\ell$ and satisfies the zero mean condition. We will denote
$\left(D\mathbf{u}\right)_{\ell}$ by $Du_{\ell}$.
We then define
\begin{align*}
\left(D^{(2)}\mathbf{u}\right)_{\ell}:=\frac{Du_{\ell}-Du_{\ell-1}}{\epsilon},
\qquad -\infty<\ell<\infty,
\end{align*}
and we define $\left(D^{(3)}\mathbf{u}\right)_{\ell}$ and
$\left(D^{(4)}\mathbf{u}\right)_{\ell}$ in a similar way. To make
the formulas concise and more readable, we sometimes denote $Du_{\ell}$ by $u'_{\ell}$,
$D^{(2)}u_{\ell}$ by $u''_{\ell}$, etc., when there is no confusion
in the expressions.

For a displacement $\mathbf{u}\in \mathcal{U}$ and its discrete derivatives, we define the
discrete $\ell_{\epsilon}^{2}$ norms by
\begin{align*}
\|\mathbf{u}\|_{\ell_{\epsilon}^{2}}&:= \left( \epsilon
\sum_{\ell=-N+1}^{N}|u_{\ell}|^{2}\right)^{1/2},\qquad
\|\mathbf{u}'\|_{\ell_{\epsilon}^{2}}:= \left( \epsilon
\sum_{\ell=-N+1}^{N}|u_{\ell}'|^{2}\right)^{1/2},\text{ etc.}
\end{align*}
Finally, for smooth real-valued functions $\E(\yb)$ defined for
$\yb\in\mathcal{Y}_{F},$ we define the first and second derivatives (variations) by
\begin{equation*}
\begin{split}
\la\delta\mathcal{E}(\mathbf{y}),\mathbf{w}\ra&:=\sum_{\ell=-N+1}^{N}
 \frac{\partial \mathcal{E}}{\partial y_{\ell}}(\mathbf{y})w_{\ell}\qquad
\text{for all }\mathbf{w}\in \mathcal{U}\\
\la\delta^2\mathcal{E}(\mathbf{y})\mathbf{v},\mathbf{w}\ra&:=\sum_{\ell,\,m=-N+1}^{N}
 \frac{\partial^2 \mathcal{E}}{\partial y_{\ell}\partial y_{m}}(\mathbf{y})v_{\ell}
 w_{m}\qquad
\text{for all }\mathbf{v},\,\mathbf{w}\in \mathcal{U}.
\end{split}
\end{equation*}

\section{The Embedded Atom Model and Its QNL Approximation}\label{qnl}
We first give a description of the next-nearest neighbor EAM Model.
\subsection{The Next-Nearest-Neighbor Embedded Atom Model}
The total energy per period of the next-nearest neighbor EAM model is
\begin{equation}\label{AtomTotalEnergy}
 \mathcal{E}_{tot}^{a}(\mathbf{y}):=\mathcal{E}^{a}(\mathbf{y})+\mathcal{F}(\mathbf{y})
\end{equation}
for deformations $\mathbf{y}\in \mathcal{Y}_{F}$
where $\mathcal{E}^{a}(\mathbf{y})$ is the
total atomistic energy and $\mathcal{F}(\mathbf{y})$
is the total external potential energy.
The total atomistic energy is the sum of the {\it embedding energy}, $\hat{\mathcal{E}}^{a}(\mathbf{y}),$
and the {\it pair potential energy}, $\tilde{\mathcal{E}}^{a}(\mathbf{y}):$
\begin{equation}\label{AtomModel}
\mathcal{E}^{a}(\mathbf{y}):=\hat{\mathcal{E}}^{a}(\mathbf{y})+\tilde{\mathcal{E}}^{a}(\mathbf{y}).
\end{equation}
The embedding energy is
\[
\hat{\mathcal{E}}^{a}(\mathbf{y}):= \epsilon \sum_{\ell=-N+1}^{N}
                            G\left(\bar\rho^a_\ell(\mathbf{y})\right)
\]
where $G(\bar\rho)$ is the embedding energy function, the total electron density
$\bar\rho^a_\ell(\mathbf{y})$ at
atom $\ell$ is
\[
\bar\rho^a_\ell(\mathbf{y}):=
\rho(y'_{\ell})+\rho(y'_{\ell}+y'_{\ell-1})
                            +\rho(y'_{\ell+1})+\rho(y'_{\ell+1}+y'_{\ell+2}),
\]
and $\rho(r/\eps)$ is the electron density contributed by an atom at distance $r.$
The pair potential energy is
\[
 \tilde{\mathcal{E}}^{a}(\mathbf{y}):=\epsilon \sum_{\ell=-N+1}^{N}
 \half \left[\phi(y'_{\ell})+\phi(y'_{\ell}+y'_{\ell-1}) +\phi(y'_{\ell+1})+\phi(y'_{\ell+1}+y'_{\ell+2})\right]
\]
where $\eps\phi(r/\eps)$ is
the pair potential interaction energy~\cite{Foiles}.
Our formulation
allows general nonlinear external potential energies $\mathcal{F}(\mathbf{y})$ defined
for $\mathbf{y}\in \mathcal{Y}_{F},$ but we note that the
total external potential
energy for periodic dead loads $\mathbf{f}$ is given by
\[
 \mathcal{F}(\mathbf{y}):=-\sum_{\ell=-N+1}^{N}\epsilon f_{\ell}y_{\ell}.
\]

The equilibrium solution $\mathbf{y}^{a}$ of the
EAM atomistic model \eqref{AtomTotalEnergy} then satisfies
\begin{align}\label{AtomEqulibriumEq2}
- \la\delta \mathcal{E}^{a}(\mathbf{y}^a),\mathbf{w}\ra=
-\la\delta\hat{\mathcal{E}}^{a}(\mathbf{y}^a),\mathbf{w}\ra
-\la\delta\tilde{\mathcal{E}}^{a}(\mathbf{y}^a),\mathbf{w}\ra
=\la\delta\mathcal{F}(\mathbf{y}^{a}),\mathbf{w}\ra \qquad \text{for all }
\mathbf{w}\in \mathcal{U}.
 \end{align}
Here the negative of the embedding force of \eqref{AtomEqulibriumEq2} is given by
\begin{align*}
\la\delta\hat{\mathcal{E}}^{a}(\mathbf{y}^a),\mathbf{w}\ra&=\epsilon
\sum_{\ell=-N+1}^{N}G'\Big(\bar\rho^a_\ell(\mathbf{y}^a) \Big)\cdot\Big[
\rho'(Dy^{a}_{\ell})w'_{\ell}+\rho'(Dy^{a}_{\ell}+Dy^{a}_{\ell-1})(w'_{\ell}+w'_{\ell-1})\Big.\\
&\qquad\qquad\qquad
\big.+\rho'(Dy^{a}_{\ell+1})w'_{\ell+1}+\rho'(Dy^{a}_{\ell+1}+Dy^{a}_{\ell+2})(w'_{\ell+1}+w'_{\ell+2})
\Big],
\end{align*}
the negative of the pair potential force of \eqref{AtomEqulibriumEq2} is given by
\begin{align*}
\la\delta\tilde{\mathcal{E}}^{a}(\mathbf{y}^a),\mathbf{w}\ra&=
\epsilon\sum_{\ell=-N+1}^{N}\half
\Big[
\phi'(Dy^{a}_{\ell})w'_{\ell}+\phi'(Dy^{a}_{\ell}+Dy^{a}_{\ell-1})(w'_{\ell}+w'_{\ell-1})
\\
&\qquad\qquad\qquad\qquad\qquad
+\phi'(Dy^{a}_{\ell+1})w'_{\ell+1}+\phi'(Dy^{a}_{\ell+1}+Dy^{a}_{\ell+2})(w'_{\ell+1}+w'_{\ell+2})
\Big]
\end{align*}
and the external force is given by
\[
 \la\delta\mathcal{F}(\mathbf{y}),\mathbf{w}\ra=\sum_{\ell=-N+1}^{N}
 \frac{\partial \mathcal{F}}{\partial y_{\ell}}(\mathbf{y})w_{\ell}\qquad
\text{for all }\mathbf{w}\in \mathcal{U}.
\]


\subsection{The Quasi-Nonlocal EAM Approximation
for Next-Nearest-Neighbor Interactions}
Hybrid atomistic-to-continuum
methods can give an accurate and efficient solution
if the deformation $\mathbf{y}\in\mathcal{Y}_{F}$ is
''smooth'' in most of the computational domain, but not in the
remaining domain where defects occur~\cite{ortner:qnl1d,dobs-qcf2}.
The goal of QC methods is to decompose the
reference lattice into an atomistic region with defects and a
continuum region with long-range elastic effects. It applies an
atomistic model to the atomistic region for accuracy and
a continuum model to the continuum region for efficiency.

In this paper, we will consider an atomistic region defined by the
atoms with reference positions $x_{\ell}$ for
$\ell=-K,\dots,K$, and a continuum region for $\ell\in
\{-N+1,\dots,-(K+3)\}\cup\{(K+3),\dots,N\}$. To eliminate the ghost
force that energy-based quasicontinuum approximations can have
\cite{Shimokawa:2004,Miller:2003a,
Ortiz:1995a,Dobson:2008b},
we define the remaining atoms, $\pm (K+1), \pm (K+2)$, to be
quasi-nonlocal atoms \cite{Shimokawa:2004,Dobson:2008b}.
For the pair potential energy, the quasi-nonlocal atoms $\pm (K+1),\pm(K+2)$ interact without
approximation with atoms in the atomistic region,
but interact through the continuum Cauchy-Born approximation
with all other
atoms~\cite{Shimokawa:2004}.  The interactions of the quasi-nonlocal atoms
for the embedding energy is slightly more complex, as given in \cite{Shimokawa:2004}
and below.

The atomistic energy associated with each atom is given by
\begin{align*}
 \mathcal{E}^{a}_{\ell}(\mathbf{y}):=
 \hat{\mathcal{E}}^{a}_{\ell}(\mathbf{y})+\tilde{\mathcal{E}}^{a}_{\ell}(\mathbf{y})
 =
 G\left(\rhoa\right)+\half \left[\phi(y'_{\ell})+\phi(y'_{\ell}+y'_{\ell-1})
 +\phi(y'_{\ell+1})+\phi(y'_{\ell+1}+y'_{\ell+2})\right]
\end{align*}
where $\hat{\mathcal{E}}^{a}_{\ell}(\mathbf{y})$ denotes
the embedding energy at atom $\ell$ and
$\tilde{\mathcal{E}}^{a}_{\ell}(\mathbf{y})$ denotes
the pair potential energy at atom $\ell$
($\hat{\mathcal{E}}^{c}_{\ell}(\mathbf{y}),$
$\hat{\mathcal{E}}^{qnl}_{\ell}(\mathbf{y}),$
$\tilde{\mathcal{E}}^{c}_{\ell}(\mathbf{y})$
and $\tilde{\mathcal{E}}^{qnl}_{\ell}(\mathbf{y})$ will
be defined analogously below),
and
the continuum energy associated with each
atom is given by
\begin{align*}
 \mathcal{E}^{c}_{\ell}(\mathbf{y}):=
 \hat{\mathcal{E}}^{c}_{\ell}(\mathbf{y})+\tilde{\mathcal{E}}^{c}_{\ell}(\mathbf{y})
=&
\half G\left(\rhoc)\right)+\half
G\left(\rhocc\right) \\
&\quad
+\half \left[\phi(y'_{\ell})+\phi(2y'_{\ell})
+\phi(y'_{\ell+1})+\phi(2y'_{\ell+1})\right]
\end{align*}
where the total continuum electron density at atom $\ell$ is
\[
\rhoc:=2\rho(y'_{\ell})
 +2\rho(2y'_{\ell}).
\]

To define the QNL energy for the quasi-nonlocal atoms, we define the
QNL electron density at atom $\ell$ by
\[
\rhoqnl:=2\rho(y'_{\ell})+2\rho(y'_{\ell}+y'_{\ell-1}).
\]
We then define the QNL energy for the quasi-nonlocal atoms by
\begin{align*}
 \mathcal{E}^{qnl}_{K+1}(\mathbf{y}):&=\hat{\mathcal{E}}^{qnl}_{K+1}
(\mathbf{y})+\tilde{\mathcal{E}}^{qnl}_{K+1}(\mathbf{y})\\
&= \half G\Big(\bar\rho^{qnl}_{K+1}(\mathbf{y})\Big)+\half G\Big(\bar\rho^c_{K+2}(\mathbf{y})\Big)\\
&\qquad+ \half\left[\phi(y'_{K+1})+\phi(y'_{K+2})+\phi(y'_{K+1}+y'_{K})+
 \phi(2y'_{K+2})\right]
\end{align*}
and
\begin{align*}
 \mathcal{E}^{qnl}_{K+2}(\mathbf{y}):&=\hat{\mathcal{E}}^{qnl}_{K+2}(\mathbf{y})
+\tilde{\mathcal{E}}^{qnl}_{K+2}(\mathbf{y})\\
&=  \half G\Big(\bar\rho^{qnl}_{K+2}(\mathbf{y})\Big)+\half G\Big(\bar\rho^c_{K+3}(\mathbf{y})\Big)\\
&\qquad+ \half\left[\phi(y'_{K+2})+\phi(y'_{K+3})+\phi(y'_{K+2}+y'_{K+1})+\phi(2y'_{K+3})\right].
\end{align*}
We define the QNL energy in a symmetric way and so only
give the formulas for $0\le \ell \le N.$

The total energy per period of the QNL model
is then given by
\begin{equation}\label{QNLTotalEnergy}
\begin{split}
 \mathcal{E}_{tot}^{qnl}(\mathbf{y})&
 :=\epsilon\sum_{\ell=-N+1}^{N}\mathcal{E}^{qnl}_{\ell}(\mathbf{y})+\mathcal{F}(\mathbf{y})\\
 &=\mathcal{E}^{qnl}(\mathbf{y})+\mathcal{F}(\mathbf{y})
= \hat{\mathcal{E}}^{qnl}(\mathbf{y})+\tilde{\mathcal{E}}^{qnl}(\mathbf{y})+\mathcal{F}(\mathbf{y}),
\end{split}
\end{equation}
where
\begin{align*}
 \mathcal{E}^{qnl}_{\ell}(\mathbf{y}):=\begin{cases}
                           &\mathcal{E}^{a}_{\ell}(\mathbf{y})\quad \text{for}\quad 0\le\ell<K+1, \\
 & \mathcal{E}^{qnl}_{\ell}(\mathbf{y})
                          \quad \text{for}\quad \ell=K+1,\,K+2,\\
&\mathcal{E}^{c}_{\ell}(\mathbf{y})
                          \quad \text{for}\quad K+2<\ell<N. \\
                          \end{cases}
\end{align*}
The equilibrium solution $\mathbf{y}^{qnl}$ of the
EAM-QNL model \eqref{QNLTotalEnergy} then satisfies
\begin{align}\label{QNLEqulibriumEq2}
-\la\delta \mathcal{E}^{qnl}(\mathbf{y}^{qnl}),\mathbf{w}\ra=
-\la\delta\hat{\mathcal{E}}^{qnl}(\mathbf{y}^{qnl}),\mathbf{w}\ra
-\la\delta\tilde{\mathcal{E}}^{qnl}(\mathbf{y}^{qnl}),\mathbf{w}\ra
=\la\delta \F(\mathbf{y}^{qnl}),\mathbf{w}\ra \qquad \text{for all }
\mathbf{w}\in \mathcal{U},
 \end{align}
 where the negative of the embedding force is given by
\begin{align}\label{QNLvariationFirst}
\begin{split}
&\la\delta \hat{\mathcal{E}}^{qnl}(\mathbf{y}^{qnl}),\mathbf{w}\ra
 =\dots\\
 &\qquad+\epsilon\sum_{\ell=0}^{K}G'\big(\bar\rho^a_\ell(\mathbf{y}^{qnl}) \big)\cdot\left[
\rho'(Dy^{qnl}_{\ell})w'_{\ell}+\rho'(Dy^{qnl}_{\ell}+Dy^{qnl}_{\ell-1})(w'_{\ell}+w'_{\ell-1})\right.\\
&\qquad\qquad\qquad\qquad\qquad\qquad\left.+\rho'(Dy^{qnl}_{\ell+1})w'_{\ell+1}
+\rho'(Dy^{qnl}_{\ell+1}+Dy^{qnl}_{\ell+2})(w'_{\ell+1}+w'_{\ell+2})
\right]\\
&\qquad+\epsilon G'\big(\bar\rho^{qnl}_{K+1}(\mathbf{y}^{qnl})\big)\cdot\left[
\rho'(Dy^{qnl}_{K+1})w'_{K+1}+\rho'(Dy^{qnl}_{K+1}+Dy^{qnl}_{K})(w'_{K+1}+w'_{K})\right]\\
&\qquad+\epsilon G'\big(\bar\rho^c_{K+2}(\mathbf{y}^{qnl})\big)
\cdot\left[
\rho'(Dy^{qnl}_{K+2})w'_{K+2}+2\rho'(2Dy^{qnl}_{K+2})(w'_{K+2})\right]\\
&\qquad +\epsilon G'\big(\bar\rho^{qnl}_{K+2}(\mathbf{y}^{qnl})\big)
\cdot\left[
\rho'(Dy^{qnl}_{K+2})w'_{K+2}+\rho'(Dy^{qnl}_{K+2}+Dy^{qnl}_{K+1})(w'_{K+2}+w'_{K+1})\right]\\
&\qquad +\epsilon G'\big(\bar\rho^c_{K+3}(\mathbf{y}^{qnl})\big)
\cdot\left[
\rho'(Dy^{qnl}_{K+3})w'_{K+3}+2\rho'(2Dy^{qnl}_{K+3})(w'_{K+3})\right]\\
&\qquad+\epsilon \sum_{\ell=K+3}^{N}\left\{
G'\big(\bar\rho^c_{\ell}(\mathbf{y}^{qnl})\big) \cdot\left[
\rho'(Dy^{qnl}_{\ell})w'_{\ell}+2\rho'(2Dy^{qnl}_{\ell})(w'_{\ell})\right]\right.\\
&\qquad\qquad\qquad\qquad\left.
+G'\big(\bar\rho^c_{\ell+1}(\mathbf{y}^{qnl})\big)
\cdot\left[
\rho'(Dy^{qnl}_{\ell+1})w'_{\ell+1}+2\rho'(2Dy^{qnl}_{\ell+1})(w'_{\ell+1})\right]\right\},
\end{split}
\end{align}
and the negative of the pair potential force is given by
\begin{align}\label{QNLvariationSecond}
\begin{split}
& \la\delta\tilde{\mathcal{E}}^{qnl}(\mathbf{y}^{qnl}),\mathbf{w}\ra
=\dots\\
&\qquad+\epsilon\sum_{\ell=0}^{K}\half\left[
\phi'(Dy^{qnl}_{\ell})w'_{\ell}+\phi'(Dy^{qnl}_{\ell}+Dy^{qnl}_{\ell-1})(w'_{\ell}+w'_{\ell-1})\right.\\
&\qquad\qquad\qquad\qquad\qquad\left.+\phi'(Dy^{qnl}_{\ell+1})w'_{\ell+1}+\phi'(Dy^{qnl}_{\ell+1}+Dy^{qnl}_{\ell+2})(w'_{\ell+1}+w'_{\ell+2})
\right]\\
&\qquad +\frac{\epsilon}{2}\left[
\phi'(Dy^{qnl}_{K+1})w'_{K+1}+\phi'(Dy^{qnl}_{K+1}+Dy^{qnl}_{K})(w'_{K+1}+w'_{K})\right]\\
&\qquad\qquad+\frac{\epsilon}{2} \left[
\phi'(Dy^{qnl}_{K+2})w'_{K+2}+2\phi'(2Dy^{qnl}_{K+2})(w'_{K+2})\right]\\
&\qquad +\frac{\epsilon}{2} \left[
\phi'(Dy^{qnl}_{K+2})w'_{K+2}+\phi'(Dy^{qnl}_{K+2}+Dy^{qnl}_{K+1})(w'_{K+2}+w'_{K+1})\right]\\
&\qquad\qquad +\frac{\epsilon}{2} \left[
\phi'(Dy^{qnl}_{K+3})w'_{K+3}+2\phi'(2Dy^{qnl}_{K+3})(w'_{K+3})\right]\\
&\qquad+\epsilon \sum_{\ell=K+3}^{N}\half\left[
\phi'(Dy^{qnl}_{\ell})w'_{\ell}
+2\phi'(2Dy^{qnl}_{\ell})w'_{\ell}
 +\phi'(Dy^{qnl}_{\ell+1})w'_{\ell+1}+2\phi'(2Dy^{qnl}_{\ell+1})w'_{\ell+1}\right].
 \end{split}
\end{align}
\section{Stability Analysis of The Atomistic and EAM-QNL Models}\label{sharp}
In this section, we will give a stability analysis for the
atomistic model and the EAM-QNL model for the next-nearest neighbor case. We will use
techniques similar to those presented in \cite{doblusort:qce.stab}
for the atomistic and QNL method for pair potentials.
\subsection{The Atomistic Model.}
The uniform deformation $\mathbf{y}_{F}$ is an equilibrium of the
atomistic model~\eqref{AtomModel}, therefore, we say that the equilibrium
$\mathbf{y}_{F}$ is stable in the atomistic model if and only if
$\la\delta^{2}\mathcal{E}^{a}(\mathbf{y}_{F})$ is positive definite,
that is,
\begin{align}\label{AtomStabEq0}
 \la\delta^{2}\mathcal{E}^{a}(\mathbf{y}_{F})\mathbf{u},\mathbf{u}\ra=
\la\delta^{2}\hat{\mathcal{E}}^{a}(\mathbf{y}_{F})\mathbf{u},\mathbf{u}\ra
+\la\delta^{2}\tilde{\mathcal{E}}^{a}(\mathbf{y}_{F})\mathbf{u},\mathbf{u}\ra>0
\quad \text{for all } \mathbf{u}\in \mathcal{U}\setminus \{\mathbf{0}\}.
\end{align}
Note that
$\la\delta^{2}\tilde{\mathcal{E}}_{a}(\mathbf{y}_{F})\mathbf{u},\mathbf{u}\ra$
is given by formula $(7)$ in \cite{doblusort:qce.stab}:
\begin{align}\label{PairAtomStabEq1}
 \begin{split}
\la\delta^{2}\tilde{\mathcal{E}}^{a}(\mathbf{y}_{F})\mathbf{u},\mathbf{u}\ra
=\tilde{A}_{F}\|D\mathbf{u}\|_{\ell_{\epsilon}^2}^2
-\epsilon^2\phi''_{2F}\|D^{(2)}\mathbf{u}\|_{\ell_{\epsilon}^2}^2,
 \end{split}
\end{align}
where
\begin{equation}\label{PairAtomCond1}
\tilde{A}_{F}:=\phi''_{F}+4\phi''_{2F}\quad\text{for}\quad
\phi''_{F}:=\phi''(F)\text{ and } \phi''_{2F}:=\phi''(2F)
\end{equation}
is the {\it continuum elastic modulus for the pair interaction potential}.
 Thus,
we only need to focus on
$\la\delta^{2}\hat{\mathcal{E}}^{a}(\mathbf{y}_{F})\mathbf{u},\mathbf{u}\ra$,
that is,
\begin{align}\label{AtomStabEq1}
\begin{split}
 \la\delta^{2}\hat{\mathcal{E}}^{a}(\mathbf{y}_{F})\mathbf{u},\mathbf{u}\ra
&=\epsilon \sum_{\ell=-N+1}^{N}\bigg\{G''_{F}\,
                          \left[\rho'_{F}(u'_{\ell}+u'_{\ell+1})
                         +\rho'_{2F}(u'_{\ell-1}+u'_{\ell}+u'_{\ell+1}+u'_{\ell+2})\right]^2\\
&\qquad\qquad\qquad\quad\left. +G'_{F}\left[\rho''_{F}(u'_{\ell})^2+\rho''_{2F}(u'_{\ell}+u'_{\ell-1})^2
              +\rho''_{F}(u'_{\ell+1})^2\right.\right.\\
&\qquad\qquad\qquad\qquad\qquad\quad\left.+\rho''_{2F}(u'_{\ell+1}+u'_{\ell+2})^2\right]\bigg\},
\end{split}
\end{align}
where
\begin{gather*}
\rho'_{F}:=\rho'(F),\quad \rho''_{F}:=\rho''(F),
\quad \rho'_{2F}:=\rho(2F),\quad \rho''_{2F}:=\rho''(2F),\\
 G'_{F}:=G'(\bar\rho^a_\ell(\mathbf{y}_F))=G'(\bar\rho^c_\ell(\mathbf{y}_F))
 =G'(\bar\rho^{qnl}_\ell(\mathbf{y}_F)),\\
 G''_{F}:=G''(\bar\rho^a_\ell(\mathbf{y}_F))
 =G''(\bar\rho^c_\ell(\mathbf{y}_F))
 =G''(\bar\rho^{qnl}_\ell(\mathbf{y}_F)).
\end{gather*}

We calculate the identities
\begin{align}\label{DerivativeEquivalence}
 \left(u'_{\ell}+u'_{\ell+1}\right)^2 &= 2\left(u'_{\ell}\right)^2+2\left(u'_{\ell+1}\right)^2-\epsilon^2(u''_{\ell+1})^2,\\
\left( u'_{\ell}+u'_{\ell+1}+u'_{\ell+2} \right)^2&= 3\left(u'_{\ell}\right)^2+3\left(u'_{\ell+1}\right)^2
                    + 3\left(u'_{\ell+2}\right)^2-3\epsilon^{2}\left(u''_{\ell+1}\right)^2
                  -3\epsilon^2\left(u''_{\ell+2}\right)^2+\epsilon^4\left(u^{(3)}_{\ell+2}\right)^2.\notag\\
2\left(u'_{\ell}+u'_{\ell+1}\right)\cdot & \left(u'_{\ell-1}
+u'_{\ell}+u'_{\ell+1}+u'_{\ell+2}\right) \notag\\
&\,=2\left[\left(u'_{\ell-1}\right)^2+3\left(u'_{\ell}\right)^2+3\left(u'_{\ell+1}\right)^2+\left(u'_{\ell+2}\right)^2\right]\notag\\
&\qquad-3\epsilon^2\left[\left(u''_{\ell}\right)^2+2\left(u''_{\ell+1}\right)^2+\left(u''_{\ell+2}\right)^2\right]
+\epsilon^4\left[\left(u^{(3)}_{\ell+1}\right)^2+\left(u^{(3)}_{\ell+2}\right)^2\right].\notag
\end{align}
We can now calculate explicitly the first equality below and then
use \eqref{DerivativeEquivalence} (with $\mathbf{u}'$ replaced by $\mathbf{u}''$) for the second equality to obtain
\begin{align*}
\left( u'_{\ell}+u'_{\ell+1}+u'_{\ell+2} +u'_{\ell+3}\right)^2
&=4\left((u'_{\ell})^2+(u'_{\ell+1})^2+(u'_{\ell+2})^2+(u'_{\ell+3})^2\right)\\
&\quad-\epsilon^2\left(u''_{\ell+1}\right)^2-\epsilon^2\left(u''_{\ell+2}\right)^2-\epsilon^2\left(u''_{\ell+3}\right)^2
-\epsilon^2\left(u''_{\ell+1}+u''_{\ell+2}\right)^2\\
&\qquad-\epsilon^2\left(u''_{\ell+2}+u''_{\ell+3}\right)^2-\epsilon^2\left(u''_{\ell+1}+u''_{\ell+2}+u''_{\ell+3}\right)^2\\
&= 4\left((u'_{\ell})^2+(u'_{\ell+1})^2+(u'_{\ell+2})^2+(u'_{\ell+3})^2\right)\\
&\quad-\epsilon^2 \left(6(u''_{\ell+1})^2+8(u''_{\ell+2})^2+6(u''_{\ell+3})^2\right)\\
&\qquad+\epsilon^4 \left(4(u^{(3)}_{\ell+2})^2+4(u^{(3)}_{\ell+3})^2\right)-\epsilon^6(u^{(4)}_{\ell+3})^2.
\end{align*}
We can then obtain from the above identities that
\begin{align}\label{AtomStabEq2}
\begin{split}
 \la\delta^{2}\hat{\mathcal{E}}^{a}(\mathbf{y}_{F})\mathbf{u},\mathbf{u}\ra&=
G''_{F}\cdot\left\{\left[4\left(\rho'_{F}\right)^2+16\left(\rho'_{2F}\right)^2
                  +16\rho'_{F}\rho'_{2F}\right]\,\|D\mathbf{u}\|_{\ell_{\epsilon}^2}^{2}\right.\\
&\qquad \qquad\left. -\epsilon^2\left[\left(\rho'_{F}\right)^2+20\left(\rho'_{2F}\right)^2
                  +12\rho'_{F}\rho'_{2F}\right]\,\|D^{(2)}\mathbf{u}\|_{\ell_{\epsilon}^2}^{2}\right.\\
&\qquad\qquad\left.+\epsilon^4\left[8\left(\rho'_{2F}\right)^2
                  +2\rho'_{F}\rho'_{2F}\right]\,\|D^{(3)}\mathbf{u}\|_{\ell_{\epsilon}^2}^{2}
                -\epsilon^6\left(\rho'_{2F}\right)^2\,\|D^{(4)}\mathbf{u}\|_{\ell_{\epsilon}^2}^{2}\right\}\\
&\qquad+G'_{F}\cdot\left\{ \left(2\rho''_{F}+8\rho''_{2F}\right)
        \,\|D\mathbf{u}\|_{\ell_{\epsilon}^2}^{2} - 2\epsilon^2\rho''_{2F}\,\|D^{(2)}\mathbf{u}\|_{\ell_{\epsilon}^2}^{2}\right\}\\
&=\left\{4G''_{F}\left(\rho'_{F}+2\rho'_{2F}\right)^2
+2G'_{F}\left(\rho''_{F}+4\rho''_{2F}\right)\right\}\|D\mathbf{u}\|_{\ell_{\epsilon}^2}^{2}\\
&\quad -\epsilon^2\left\{G''_{F}\left[\left(\rho'_{F}\right)^2
+20\left(\rho'_{2F}\right)^2+12\rho'_{F}\rho'_{2F}\right]
+G'_{F}\,2\rho''_{2F}\right\}\|D^{(2)}\mathbf{u}\|_{\ell_{\epsilon}^2}^{2}\\
&\quad+\epsilon^4 G''_{F}\left[8\left(\rho'_{2F}\right)^2+2\rho'_{F}\rho'_{2F}\right]
\|D^{(3)}\mathbf{u}\|_{\ell_{\epsilon}^2}^{2}\\
&\quad-\epsilon^6 G''_{F}\left(\rho'_{2F}\right)^2 \|D^{(4)}\mathbf{u}\|_{\ell_{\epsilon}^2}^{2}.
\end{split}
\end{align}

We define the {\it continuum elastic modulus for the embedding energy} to be
\begin{equation}\label{EAMCond1}
\hat{A}_{F}:=4G''_{F}\left(\rho'_{F}+2\rho'_{2F}\right)^2
+2G'_{F}\left(\rho''_{F}+4\rho''_{2F}\right).
\end{equation}
and
\begin{gather*}
A_{F}:=\hat{A}_{F}+\tilde{A}_{F},\quad B_{F}:= -\left[\phi''_{2F}+G''_{F}\big((\rho'_{F})^2+20(\rho'_{2F})^2+12\rho'_{F}\rho'_{2F}\big)+G'_{F}\left(2\rho''_{2F}\right)\right],\\
C_{F}:=G''_{F}\left(8(\rho'_{2F})^2+2\rho'_{F}\rho'_{2F}\right),\quad\text{and}\quad
 D_{F}:=-G''_{F}\left(\rho'_{2F}\right)^2.
\end{gather*}
Then \eqref{AtomStabEq0} becomes
\begin{equation}\label{Four}
\begin{split}
\la\delta^{2}\mathcal{E}^{a}(\mathbf{y}_{F})\mathbf{u},\mathbf{u}\ra=&
A_F\|D\mathbf{u}\|_{\ell_{\epsilon}^{2}}^{2}
+\epsilon^2B_F \|D^{(2)}\mathbf{u}\|_{\ell_{\epsilon}^{2}}^{2}
+\epsilon^4
C_F
\|D^{(3)}\mathbf{u}\|_{\ell_{\epsilon}^2}^{2}
+\epsilon^6 D_F
\|D^{(4)}\mathbf{u}\|_{\ell_{\epsilon}^2}^{2}.
\end{split}
\end{equation}

We will analyze the stability of $\la\delta^{2}\mathcal{E}^{a}(\mathbf{y}_{F})\mathbf{u},\mathbf{u}\ra$
by using the Fourier representation~\cite{HudsOrt:a}
\[
Du_{\ell}=\sum_{k=-N+1}^{N}\frac{c_{k}}{\sqrt{2}}\cdot\exp\left(i\,k\frac{\ell}{N}\pi\right).
\]
It then follows from the discrete orthogonality of the Fourier basis that
\begin{equation}\label{AtomStabFourierEq1}
\begin{split}
\la\delta^{2}\mathcal{E}^{a}(\mathbf{y}_{F})\mathbf{u},\mathbf{u}\ra &=
\sum_{k=-N+1}^{N} |c_{k}|^2\cdot \Bigg\{
A_F
+B_F
\left[4\sin^{2}\left(\frac{k\pi}{2N}\right)\right]\\
&\qquad\qquad\qquad +C_F\left[4\sin^2\left(\frac{k\pi}{2N}\right)\right]^2
+D_F\left[4\sin^2\left(\frac{k\pi}{2N}\right)\right]^3\Bigg\}.
\end{split}
\end{equation}
We then
see from \eqref{AtomStabFourierEq1} that the eigenvalues $\lambda_k$ for
$k=-N+1,\dots,N$ of
$\la\delta^{2}\mathcal{E}^{a}(\mathbf{y}_{F})\mathbf{u},\mathbf{u}\ra$
with respect to the $\|D\mathbf{u}\|_{\ell_{\epsilon}^{2}}$
norm are given by
\[
\lambda_k=\lambda_{F}(s_k)\quad\text{for }\quad s_k=4\sin^2\left(\frac{k\pi}{2N}\right)
\]
where
\[
\lambda_{F}(s):=A_{F}+B_{F}s+C_{F}s^2+D_{F}s^3.
\]

From the pair interaction potential, electron density function, and
embedding energy function given in
 Figure~2
in \cite{Foiles}, we assume that
\begin{align}\label{EAMfunAssumption1}
\begin{split}
\phi''_{F}>0,\,\phi''_{2F}<0;\quad
 \rho'_{F}\le 0, \, \rho'_{2F} \le 0;\quad
 \rho''_{F}\ge 0,\,\rho''_{2F} \ge 0;\quad\text{and}\quad
G''_{F}\ge 0.
\end{split}
\end{align}
We then have from the assumption \eqref{EAMfunAssumption1} that
\begin{align}\label{condition1}
C_{F}> 0,\quad
D_{F}< 0,\quad\text{and}\quad
8|D_{F}|\le C_{F}.\quad
\end{align}
We can check that \eqref{condition1} implies that $|D_Fs|\le 4|D_F|\le C_F/2,$ for $0\le s\le 4,$ so
\begin{equation}\label{lambda}
\lambda_{F}'(s)=B_F+2C_Fs+3D_Fs^2\ge B_{F}+\frac {C_{F}}2 s\quad\text{for all}\quad 0\le s\le 4.
\end{equation}

We conclude from \eqref{lambda} that the condition $B_F\ge 0$ or equivalently
\begin{align}\label{EAMfunAssumption2}
 \begin{split}
&\phi''_{2F}+G''_{F}\left[\left(\rho'_{F}\right)^2
+20\left(\rho'_{2F}\right)^2+12\rho'_{F}\rho'_{2F}\right]
+G'_{F}\,2\rho''_{2F}=-B_F\le 0,
 \end{split}
\end{align}
and the assumptions \eqref{EAMfunAssumption1} imply that $\lambda(s)$ is increasing for
$0\le s \le 4.$  We thus have the sharp stability result
\begin{equation}\label{eq:sharp}
\la\delta^{2}\mathcal{E}^{a}(\mathbf{y}_{F})\mathbf{u},\mathbf{u}\ra
\ge \lambda_{F}(s_1) \|D\mathbf{u}\|_{\ell_{\epsilon}^{2}}^2\ge \left(\hat{A}_{F}+\tilde{A}_{F}\right)
\|D\mathbf{u}\|_{\ell_{\epsilon}^{2}}^2\quad\text{for all }\mathbf{u}\in \mathcal{U}.
\end{equation}
We summarize this result in the following theorem:
\begin{theorem}\label{AtomStabThm}
  Suppose that the hypotheses \eqref{EAMfunAssumption1} and \eqref{EAMfunAssumption2} hold.
Then the uniform deformation $\mathbf{y}_{F}$ is stable for the
atomistic model if and only if
\begin{align*}
\lambda_{F}(s_1)&= A_F
+B_F
\left[4\sin^{2}\left(\frac{\pi}{2N}\right)\right]
+C_F\left[4\sin^2\left(\frac{\pi}{2N}\right)\right]^2
+D_F\left[4\sin^2\left(\frac{\pi}{2N}\right)\right]^3\\
  &=\hat{A}_{F}+\tilde{A}_{F}-4\sin^2\left(\frac{\pi}{2N}\right)\left\{\phi''_{2F}+G''_{F}\left[\left(\rho'_{F}\right)^2
+20\left(\rho'_{2F}\right)^2+12\rho'_{F}\rho'_{2F}\right]
+G'_{F}\,2\rho''_{2F}\right\}\\
&\qquad+4^2\sin^4\left(\frac{\pi}{2N}\right)G''_{F}
\left[\eta\left(\rho'_{2F}\right)^2+2\rho'_{F}\rho'_{2F}\right]
-4^3\sin^6\left(\frac{\pi}{2N}\right)G''_{F}\left(\rho'_{2F}\right)^2>0.
\end{align*}
\end{theorem}

\begin{remark}\label{unstab}
The role of the assumption~\eqref{EAMfunAssumption2} is to guarantee that
$u_{\ell}'=\sin(\epsilon \ell \pi)$ is the eigenfunction corresponding to
the smallest eigenvalue of $\la\delta^{2}\mathcal{E}^{a}(\mathbf{y}_{F})\mathbf{u},\mathbf{u}\ra$
with respect to the norm $\|D\mathbf{u}\|_{\ell_{\epsilon}^{2}}.$
In fact, we can see from the above Fourier analysis that
$u_{\ell}'=\sin(\epsilon \ell \pi)$ is not the
smallest eigenvalue of $\la\delta^{2}\mathcal{E}^{a}(\mathbf{y}_{F})\mathbf{u},\mathbf{u}\ra$
with respect to the norm $\|D\mathbf{u}\|_{\ell_{\epsilon}^{2}}$ for sufficiently
large N
if \eqref{EAMfunAssumption2} does not hold since then $\lambda'(0)<0.$

The assumption~\eqref{EAMfunAssumption2} on the the pair interaction potential, electron
density function, and embedding energy function cannot be expected to generally
hold for physical embedded atom models since the nearest neighbor term $G''_{F}(\rho'_{F})^2>0$
dominates. We note, however, that generally $G'_{F}<0$ for $F>1$ ~\cite{Mishin}, in which case
$G'_{F}\,2\rho''_{2F}<0;$ so \eqref{EAMfunAssumption2} is more likely to hold for tensile
strains $F>1.$
\end{remark}
\subsection{The EAM-QNL Model.}
Now we will analyze the stability of the EAM-QNL model for next-nearest neighbor interactions.
The Fourier techniques used to analyze the stability of the atomistic model
cannot be used for the EAM-QNL model
because the Fourier modes are no longer eigenfunctions.
Recall that the total atomistic interaction energy of
the QNL model is
$
 \mathcal{E}^{qnl}(\mathbf{y})
:=\hat{\mathcal{E}}^{qnl}(\mathbf{y})+\tilde{\mathcal{E}}^{qnl}(\mathbf{y})=\epsilon\sum_{\ell=-N+1}^{N}\mathcal{E}^{qnl}_{\ell}(\mathbf{y}),
$ where $\mathcal{E}^{qnl}_{\ell}(\mathbf{y})$ is symmetric in
$\ell\in \{-N+1,\dots,N\}$ and is given by
\begin{align*}
 \mathcal{E}^{qnl}_{\ell}(\mathbf{y}):=\begin{cases}
                           &\mathcal{E}^{a}_{\ell}(\mathbf{y})\quad \text{for}\quad 0\le\ell<K+1, \\
 & \mathcal{E}^{qnl}_{K+1}(\mathbf{y})
                          \quad \text{for}\quad \ell=K+1,\\
& \mathcal{E}^{qnl}_{K+2}(\mathbf{y})
                          \quad \text{for}\quad \ell=K+2,\\
&\mathcal{E}^{c}_{\ell}(\mathbf{y})
                          \quad \text{for}\quad K+2<\ell<N. \\
                          \end{cases}
\end{align*}
Since the QNL energy is consistent (see the consistency error analysis
in Section~\ref{convergence}), $\mathbf{y}_{F}$ is still an
equilibrium of $\mathcal{E}^{qnl}(\mathbf{y})$ \cite{Shimokawa:2004}. Therefore,
we will focus on
$\la\delta^{2}\mathcal{E}^{qnl}(\mathbf{y}_{F})\mathbf{u},\mathbf{u}\ra$ to estimate the
stability. The second variation of $\mathcal{E}^{qnl}(\mathbf{y})$
evaluated at $\mathbf{y}=\mathbf{y}_{F}$ is given by
\begin{equation}\label{QNLStabEq0}
{\la\delta^{2}\mathcal{E}}^{qnl}(\mathbf{y}_{F})\mathbf{u},\mathbf{u}\ra
= \la\delta^{2}\hat{\mathcal{E}}^{qnl}(\mathbf{y}_{F})\mathbf{u},\mathbf{u}\ra
+ \la\delta^{2}\tilde{\mathcal{E}}^{qnl}(\mathbf{y}_{F})\mathbf{u},\mathbf{u}\ra.
\end{equation}
We first compute the second term of \eqref{QNLStabEq0} and get
\begin{align}\label{QNLStabPair1}
 \begin{split}
\la\delta^{2}\tilde{\mathcal{E}}^{qnl}&(\mathbf{y}_{F})\mathbf{u},\mathbf{u}\ra\\
&=\epsilon\sum_{\ell=-K}^{K}\half\left\{\phi''_{F}\left[\left(u'_{\ell}\right)^2+\left(u'_{\ell+1}\right)^2\right]
+\phi''_{2F}\left[\left(u'_{\ell}+u'_{\ell-1}\right)^2+\left(u'_{\ell+1}+u'_{\ell+2}\right)^2\right]\right\}\\
&\qquad\quad+\frac{\epsilon}{2}\left\{\phi''_{F}\left[\left(u'_{K+1}\right)^2+\left(u'_{K+2}\right)^2\right]
                     +\phi''_{2F}\left[\left(u'_{K+1}+u'_{K}\right)^2+4\left(u'_{K+2}\right)^2\right]\right\}\\
&\qquad\quad+\frac{\epsilon}{2}\left\{\phi''_{F}\left[\left(u'_{K+2}\right)^2+\left(u'_{K+3}\right)^2\right]
                     +\phi''_{2F}\left[\left(u'_{K+2}+u'_{K+1}\right)^2+4\left(u'_{K+3}\right)^2\right]\right\}\\
&\quad+\dots+\epsilon\sum_{\ell=K+3}^{N}\half\left\{\phi''_{F}\left[\left(u'_{\ell}\right)^2+\left(u'_{\ell+1}\right)^2\right]
+\phi''_{2F}\left[4\left(u'_{\ell} \right)^2+4\left(u'_{\ell+1} \right)^2\right]\right\}.
 \end{split}
\end{align}
Here we omit the terms whose indices $\ell\in\{-N+1,\dots,-(K+3)\}$
since the QNL energy is symmetric.  Then we compute the first term,
which is given by the following expression:
\begin{equation}\label{QNLStabEq1}
 \begin{split}
 &\la\delta^{2}\hat{\mathcal{E}}^{qnl}(\mathbf{y}_{F})\mathbf{u},\mathbf{u}\ra
=\dots \\
&\quad+\epsilon\sum_{\ell=0}^{K}
\Big\{G''_{F}\,\left[\rho'_{F}(u'_{\ell}+u'_{\ell+1})
                         +\rho'_{2F}(u'_{\ell-1}+u'_{\ell}+u'_{\ell+1}+u'_{\ell+2})\right]^2\\
&\quad +G'_{F}\left[\rho''_{F}(u'_{\ell})^2+\rho''_{2F}(u'_{\ell}+u'_{\ell-1})^2
              +\rho''_{F}(u'_{\ell+1})^2+\rho''_{2F}(u'_{\ell+1}+u'_{\ell+2})^2\right]\Big\}\\
&\quad +2\epsilon
G''_{F}\left[\rho'_{F}u'_{K+1}+\rho'_{2F}\left(u'_{K+1}+u'_{K}\right)\right]^2
+\epsilon G'_{F}\left[\rho''_{F}(u'_{K+1})^2+\rho''_{2F}\left(u'_{K+1}+u'_{K}\right)^2\right]\\
&\quad+2\epsilon
G''_{F}\left(\rho'_{F}+2\rho'_{2F}\right)^2(u'_{K+2})^2
+\epsilon G'_{F}\left(\rho''_{F}+4\rho''_{2F}\right)(u'_{K+2})^2\\
&\quad +2\epsilon
G''_{F}\left[\rho'_{F}u'_{K+2}+\rho'_{2F}\left(u'_{K+2}+u'_{K+1}\right)\right]^2
+\epsilon G'_{F}\left[\rho''_{F}(u'_{K+2})^2+\rho''_{2F}\left(u'_{K+2}+u'_{K+1}\right)^2\right]\\
&\quad+2\epsilon
G''_{F}\left(\rho'_{F}+2\rho'_{2F}\right)^2(u'_{K+3})^2
+\epsilon G'_{F}\left(\rho''_{F}+4\rho''_{2F}\right)(u'_{K+3})^2\\
&\quad+\epsilon
\sum_{\ell=K+3}^{N}\left[2G''_{F}\left(\rho'_{F}+2\rho'_{2F}\right)^2
+G'_{F}\left(\rho''_{F}+4\rho''_{2F}\right)\right]\left[(u'_{\ell})^2+(u'_{\ell+1})^2\right].
\end{split}
\end{equation}
Now we use \eqref{DerivativeEquivalence} again to rewrite
\eqref{QNLStabEq1} in the following form
\begin{align*}
 \la\delta^2\hat{\mathcal{E}}^{qnl}(\mathbf{y}_{F})\mathbf{u},\mathbf{u}\ra&=
\epsilon
\sum_{\ell=-N+1}^{N}\left[2G''_{F}\left(\rho'_{F}+2\rho'_{2F}\right)^2
+G'_{F}\left(\rho''_{F}+4\rho''_{2F}\right)\right]\left[(u'_{\ell})^2+(u'_{\ell+1})^2\right]\\
&\quad +\dots -\epsilon^3\sum_{\ell=0}^{K}\left\{G''_{F}\cdot
\left[(\rho'_{F})^2+20(\rho'_{2F})^2+12\rho'_{F}\rho'_{2F}\right]+G'_{F}\cdot 2\rho''_{2F}\right\}\left(D^{(2)}u_{\ell}\right)^2\\
&\quad -\epsilon^3 \left\{G''_{F}\cdot
\left[(\rho'_{F})^2+16(\rho'_{2F})^2+11\rho'_{F}\rho'_{2F}\right]+G'_{F}\cdot 2\rho''_{2F}\right\}\left(D^{(2)}u_{K+1}\right)^2\\
&\quad -\epsilon^3 \big\{G''_{F}\cdot
\left[8(\rho'_{2F})^2+5\rho'_{F}\rho'_{2F}\right]+G'_{F}\cdot 2\rho''_{2F}\big\}\left(D^{(2)}u_{K+2}\right)^2\\
&\quad+\epsilon^5\sum_{\ell=0}^{K+1}G''_{F}\cdot\left[8(\rho'_{2F})^2+2\rho'_{F}\rho'_{2F}\right]\left(D^{(3)}u_{\ell}\right)^2\\
&\quad+\epsilon^5G''_{F}\cdot\left[4(\rho'_{2F})^2+\rho'_{F}\rho'_{2F}\right]\left(D^{(3)}u_{K+2}\right)^2
-\epsilon^7\sum_{\ell=0}^{K+2}G''_{F}\cdot(\rho'_{2F})^2\left(D^{(4)}u_{\ell}\right)^2.
\end{align*}
Combining $ \la\delta^2\hat{\mathcal{E}}^{qnl}(\mathbf{y}_{F})\mathbf{u},\mathbf{u}\ra$ and
$ \la\delta^2\tilde{\mathcal{E}}^{qnl}(\mathbf{y}_{F})\mathbf{u},\mathbf{u}\ra$ together we obtain
\begin{align*}
 \la\delta^2{\mathcal{E}}^{qnl}(\mathbf{y}_{F})\mathbf{u},\mathbf{u}\ra
&=
\epsilon\sum_{\ell=-N+1}^{N}\left(\hat{A}_{F}+\tilde{A}_{F}\right)\left(Du_{\ell}\right)^2+\dots \\
&\quad -\epsilon^3\sum_{\ell=0}^{K}\left\{\phi''_{2F}+G''_{F}\cdot
\left[(\rho'_{F})^2+20(\rho'_{2F})^2+12\rho'_{F}\rho'_{2F}\right]+G'_{F}\cdot 2\rho''_{2F}\right\}\left(D^{(2)}u_{\ell}\right)^2\\
&\quad -\epsilon^3 \left\{\phi''_{2F}+G''_{F}\cdot
\left[(\rho'_{F})^2+16(\rho'_{2F})^2+11\rho'_{F}\rho'_{2F}\right]+G'_{F}\cdot 2\rho''_{2F}\right\}\left(D^{(2)}u_{K+1}\right)^2\\
&\quad -\epsilon^3 \big\{\phi''_{2F}+G''_{F}\cdot
\left[8(\rho'_{2F})^2+5\rho'_{F}\rho'_{2F}\right]+G'_{F}\cdot 2\rho''_{2F}\big\}\left(D^{(2)}u_{K+2}\right)^2\\
&\quad+\epsilon^5\sum_{\ell=0}^{K+1}G''_{F}\cdot\left[8(\rho'_{2F})^2+2\rho'_{F}\rho'_{2F}\right]\left(D^{(3)}u_{\ell}\right)^2\\
&\quad+\epsilon^5G''_{F}\cdot\left[4(\rho'_{2F})^2+\rho'_{F}\rho'_{2F}\right]\left(D^{(3)}u_{K+2}\right)^2
-\epsilon^7\sum_{\ell=0}^{K+2}G''_{F}\cdot(\rho'_{2F})^2\left(D^{(4)}u_{\ell}\right)^2.
\end{align*}
 Because of the hypotheses \eqref{EAMfunAssumption1} and
\eqref{EAMfunAssumption2}, we have that
\begin{align*}
&\phi''_{2F}+G''_{F}\cdot\left[(\rho'_{F})^2+16(\rho'_{2F})^2+11\rho'_{F}\rho'_{2F}\right]
 +G'_{F}\cdot2\rho''_{2F}\le 0,\\
&\phi''_{2F}+G''_{F}\cdot
\left[8(\rho'_{2F})^2+5\rho'_{F}\rho'_{2F}\right]+G'_{F}\cdot2\rho''_{2F}\le
0.
\end{align*}
Thus, using
\begin{align*}
\left(D^{(4)}{u}_{\ell}\right)^2=
\left[\frac{1}{\epsilon}\left(D^{(3)}u_{\ell}-D^{(3)}{u}_{\ell-1}\right)\right]^2\le
\frac{2}{\epsilon^2}\left[\left(D^{(3)}{u}_{\ell}\right)^2+
\left(D^{(3)}{u}_{\ell-1}\right)^2\right]
\end{align*}
and noting that $G''_{F}\cdot(\rho'_{2F})^2\ge 0$,
we have
\begin{equation}\label{util}
\begin{split}
 &\epsilon^5\sum_{\ell=0}^{K+1}G''_{F}\cdot\left[8(\rho'_{2F})^2+2\rho'_{F}\rho'_{2F}\right]\left(D^{(3)}u_{\ell}\right)^2\\
&\qquad+\epsilon^5G''_{F}\cdot\left[4(\rho'_{2F})^2+\rho'_{F}\rho'_{2F}\right]\left(D^{(3)}u_{K+2}\right)^2
-\epsilon^7\sum_{\ell=0}^{K+2}G''_{F}\cdot(\rho'_{2F})^2\left(D^{(4)}u_{\ell}\right)^2\\
&\quad \ge
\epsilon^5\sum_{\ell=0}^{K+1}G''_{F}\cdot\left[4(\rho'_{2F})^2+2\rho'_{F}\rho'_{2F}\right]\left(D^{(3)}u_{\ell}\right)^2\\
&\qquad+\epsilon^5G''_{F}\cdot\left[2(\rho'_{2F})^2+\rho'_{F}\rho'_{2F}\right]\left(D^{(3)}u_{K+2}\right)^2\ge 0.
\end{split}
\end{equation}
 So, except in the case $K\in\{N-2,\dots,N\}$ when there is no continuum
region, it follows that $\mathbf{y}_{F}$ is stable in the QNL model
if and only if $\hat{A}_{F}+\tilde{A}_{F}>0$.

 Now we can give a sharp stability estimate for the QNL model from the above estimates
 and the arguments in \cite{doblusort:qce.stab,LuskinXingjie}.
\begin{theorem}\label{QNLStabThm}
 Suppose that $K<N-2$ and the hypotheses \eqref{EAMfunAssumption1} and \eqref{EAMfunAssumption2} hold,
then the uniform deformation $\mathbf{y}_{F}$ is stable in the QNL model
if and only if $\hat{A}_{F}+\tilde{A}_{F}>0$.
\end{theorem}
\begin{remark}\label{unstab2} The role of the assumption~\eqref{EAMfunAssumption2} in Theorem~\ref{QNLStabThm},
as in Theorem~\ref{AtomStabThm}, is to give a necessary condition
for $u'_{\ell}=\sin(\epsilon \ell \pi)$ to be the eigenfunction corresponding to
the smallest eigenvalue of $\la\delta^{2}\mathcal{E}^{qnl}(\mathbf{y}_{F})\mathbf{u},\mathbf{u}\ra$
with respect to the norm $\|D\mathbf{u}\|_{\ell_{\epsilon}^{2}}.$
\end{remark}

\begin{remark}
From Theorem \ref{AtomStabThm} and Theorem \ref{QNLStabThm}, we conclude that the difference
between the sharp stability conditions of the fully atomistic and QNL models is of order $O(\epsilon^2)$.
This result is the same as for the pair potential case \cite{doblusort:qcf.stab}.
\end{remark}
\begin{remark}\label{embed}
We noted in Remark~\ref{unstab} that the assumption \eqref{EAMfunAssumption2} is necessary for Theorem
\ref{AtomStabThm}. We now give an explicit example showing that
the uniform deformation can be more stable for the EAM-QCL model than for the
fully atomistic model when \eqref{EAMfunAssumption2} fails.  We recall that the EAM-QCL model is the EAM-QNL model
with no atomistic region, that is,
\begin{equation*}
 \mathcal{E}^{qcl}(\mathbf{y})
 :=\epsilon\sum_{\ell=-N+1}^{N}\mathcal{E}^{c}_{\ell}(\mathbf{y}).
\end{equation*}

We consider the case when
\begin{equation}\label{EAMfunAssumption3}
\phi''_{2F}+G''_{F}\left(\rho'_{F}+2\rho'_{2F}\right)^2+G'_{F}2\rho''_{2F}>0,
\end{equation}
which implies that \eqref{EAMfunAssumption2} does not hold since it then follows from
\eqref{EAMfunAssumption1} that
\begin{equation*}
\begin{split}
&\phi''_{2F}+G''_{F}\left[\left(\rho'_{F}\right)^2+20\left(\rho'_{2F}\right)^2
+12\rho'_{F}\rho'_{2F}\right]+G'_{F}2\rho''_{2F}\\
&\qquad=\left[\phi''_{2F}+G''_{F}\left(\rho'_{F}+2\rho'_{2F}\right)^2+G'_{F}2\rho''_{2F}\right]
+8G''_{F}\left(2\left(\rho'_{2F}\right)^2
+\rho'_{F}\rho'_{2F}\right)\\
&\qquad>0.
\end{split}
\end{equation*}

We define the oscillatory displacement $\tilde{\mathbf{u}}$ by
\[
\tilde u_{\ell}=(-1)^{\ell}\eps/(2\sqrt 2),
\]
  so
\[
\tilde u_{\ell}'=(-1)^{\ell}/(\sqrt 2),\quad \|D\tilde{\mathbf{u}}\|_{\ell^2_{\epsilon}}=1,\quad
\tilde u_{\ell}''=(-1)^{\ell}(\sqrt 2)/\eps.
\]
  We then calculate from
\eqref{PairAtomStabEq1} and \eqref{AtomStabEq1} that
\begin{equation}\label{atomstab}\begin{split}
\la\delta^2\mathcal{E}^a(\mathbf{y}_{F})\tilde{\mathbf{u}},\tilde{\mathbf{u}}\ra &=
\la\delta^2\tilde{\mathcal{E}}^a(\mathbf{y}_{F})\tilde{\mathbf{u}},\tilde{\mathbf{u}}\ra
+\la\delta^2\tilde{\mathcal{E}}^a(\mathbf{y}_{F})\tilde{\mathbf{u}},\tilde{\mathbf{u}}\ra\\
&=\epsilon\sum_{\ell=-N+1}^{N}G'_{F}2\rho''_{F}\half
+\left(\phi''_{F}+4\phi''_{2F}\right)\|D\tilde{\mathbf{u}}\|_{\ell^2_{\epsilon}}^2
+(-\epsilon^2\phi''_{2F})\|D^{(2)}\tilde{\mathbf{u}}\|^2_{\ell^2_{\epsilon}}\\
&=G'_{F}2\rho''_{F}+
\left(\phi''_{F}+4\phi''_{2F}\right)-4\phi''_{2F}=\phi''_{F}+G'_{F}2\rho''_{F}.
\end{split}
\end{equation}
Thus, we obtain that
\[
\inf_{{\mathbf{u}}\in\mathcal{U}\setminus\{\mathbf{0}\},\,\|D{\mathbf{u}}\|_{\ell^2_{\epsilon}}=1}
\la\delta^2\mathcal{E}^a(\mathbf{y}_{F}){\mathbf{u}},{\mathbf{u}}\ra\le
\phi''_{F}+G'_{F}2\rho''_{F}.
\]
On the other hand, we have that
 \begin{equation*}
\inf_{{\mathbf{u}}\in\mathcal{U}\setminus\{\mathbf{0}\},\,\|D{\mathbf{u}}\|_{\ell^2_{\epsilon}}=1}
\la\delta^2\mathcal{E}^{qcl}(\mathbf{y}_{F})\mathbf{u},\mathbf{u}\ra =\tilde{A}_{F}+\tilde{A}_{F}
=4\left[\phi''_{2F}+G''_{F}\left(\rho'_{F}+2\rho'_{2F}\right)^2+G'_{F}2\rho''_{2F}\right]
+ \phi''_{F}+G'_{F}2\rho''_{F}.
\end{equation*}
Therefore, from \eqref{EAMfunAssumption3} we have
\begin{align*}
\inf_{{\mathbf{u}}\in\mathcal{U}\setminus\{\mathbf{0}\},\,\|D{\mathbf{u}}\|_{\ell^2_{\epsilon}}=1}
\la\delta^2\mathcal{E}^{qcl}(\mathbf{y}_{F}){\mathbf{u}},{\mathbf{u}}\ra > \phi''_{F}+G'_{F}2\rho''_{F}
 \ge
\inf_{{\mathbf{u}}\in\mathcal{U}\setminus\{\mathbf{0}\},\,\|D{\mathbf{u}}\|_{\ell^2_{\epsilon}}=1}
\la\delta^2\mathcal{E}^a(\mathbf{y}_{F})\mathbf{u},\mathbf{u}\ra.
\end{align*}
This inequality indicates that the uniform deformation
$\mathbf{y}_{F}$ can be unstable for the atomistic model, but stable
for the EAM-QCL model, when the assumption \eqref{EAMfunAssumption3}
fails.

We cannot conclude from this argument, though, that the atomistic model
is less stable than the EAM-QNL model with a nontrivial atomistic region, i.e.,
$K>0.$
To see this, we consider an oscillatory displacement $\hat{\mathbf{u}}\in\mathcal{U}$
with support only in the atomistic region (a similar test function is used in \cite{belik10}):
\begin{align*}
\hat{u}_{\ell}=
\begin{cases}
&\frac{(-1)^{\ell}\epsilon}{2\sqrt{2}},\quad \ell=-(K-1),\dots,(K-1),\\
&0, \quad\text{otherwise}.
\end{cases}
\end{align*}
Then since $\hat{u}'_{\ell}=\left(\hat{u}_{\ell}-\hat{u}_{\ell-1}\right)/\epsilon$, we have
\begin{align*}
\hat{u}'_{\ell}=
\begin{cases}
&\frac{(-1)^{\ell}}{\sqrt{2}},\,\ell=-(K-2),\dots,(K-1),\\
&\frac{(-1)^{K}}{2\sqrt{2}},\,\ell=K,\\
&\frac{(-1)^{-(K-1)}}{2\sqrt{2}},\,\ell=-(K-1),\\
&0,\quad\text{otherwise}.
\end{cases}
\end{align*}
We substitute the displacement $\hat{\mathbf{u}}$ into \eqref{QNLStabEq1} and get
\begin{align}\label{RemQNLStabEq1}
\begin{split}
\la\delta^2\hat{\mathcal{E}}^{qnl}(\mathbf{y}_{F})\hat{\mathbf{u}},\hat{\mathbf{u}}\ra&
=\epsilon\sum_{\ell=-(K-2)}^{K-3}G'_{F}\rho''_{F}+2\epsilon\left\{G''_{F}\frac{1}{8}\left[3\left(\rho'_{2F}\right)^2+2\left(\rho'_{F}-\rho'_{2F}\right)^2\right]
+G'_{F}\left[\frac{7}{4}\rho''_{F}+\frac{1}{2}\rho''_{2F}\right]\right\}\\
&=\epsilon2(K-2) G'_{F}\rho''_{F}+O(\epsilon).
\end{split}
\end{align}
Similarly, we substitute $\hat{\mathbf{u}}$ into \eqref{QNLStabPair1} and get
\begin{align}\label{RemQNLStabEq0}
\begin{split}
\la\delta^2\tilde{\mathcal{E}}^{qnl}(\mathbf{y}_{F})\hat{\mathbf{u}},\hat{\mathbf{u}}\ra&
=\epsilon\sum_{\ell=-(K-2)}^{K-3}\half\phi''_{F}+O(\epsilon)=\epsilon(K-2)\phi''_{F}+O(\epsilon).
\end{split}
\end{align}
Therefore, we obtain that
\begin{align*}
\la\delta^2{\mathcal{E}}^{qnl}(\mathbf{y}_{F})\hat{\mathbf{u}},\hat{\mathbf{u}}\ra&
=\la\delta^2\hat{\mathcal{E}}^{qnl}(\mathbf{y}_{F})\hat{\mathbf{u}},\hat{\mathbf{u}}\ra
+\la\delta^2\tilde{\mathcal{E}}^{qnl}(\mathbf{y}_{F})\hat{\mathbf{u}},\hat{\mathbf{u}}\ra
=\epsilon(K-2)(\phi''_{F}+2G'_{F}\rho''_{F})+O(\epsilon).
\end{align*}
Note that
\[
\|\hat{\mathbf{u}}'\|_{\ell_{\epsilon}^2}^2=\epsilon\sum_{\ell=-N+1}^{N}(u'_{\ell})^2=\epsilon(K-1)+\frac{\epsilon}{4},
\]
Thus, we obtain from the above and \eqref{atomstab} that
\begin{align*}
\frac{\la\delta^2{\mathcal{E}}^{qnl}(\mathbf{y}_{F})\hat{\mathbf{u}},\hat{\mathbf{u}}\ra}
{\|\hat{\mathbf{u}}'\|_{\ell_{\epsilon}^2}^2}=\left(\phi''_{F}+2G'_{F}\rho''_{F}\right)
+O\left(\frac 1K\right)=
\frac{\la\delta^2{\mathcal{E}}^{a}(\mathbf{y}_{F}){\tilde{\mathbf{u}},\tilde{\mathbf{u}}\ra}}{\|\tilde{\mathbf{u}}'\|_{\ell_{\epsilon}^2}^2}
+O\left(\frac 1K\right).
\end{align*}
 This indicates that when \eqref{EAMfunAssumption3} holds and $K$ is sufficiently large, the
  EAM-QNL model is also less stable
 than the EAM-QCL model.

\end{remark}
\section{Consistency Error and Convergence of The EAM-QNL Model.}\label{convergence}

Setting
$\mathbf{y}^{qnl}=\mathbf{y}_{F}+\mathbf{u}^{qnl}$ and
$\mathbf{y}^{a}=\mathbf{y}_{F}+\mathbf{u}^{a}$, where both
$\mathbf{u}^{qnl}$ and $\mathbf{u}^{a}$ belong to $\mathcal{U},$
we define the quasicontinuum error to be
\[
\mathbf{e}^{qnl}:=\mathbf{y}^{a}-\mathbf{y}^{qnl}=\mathbf{u}^{a}-\mathbf{u}^{qnl}.
\]
To simplify the error analysis, we consider the
linearization of the atomistic equilibrium equations \eqref{AtomEqulibriumEq2} and the
associated EAM-QNL equilibrium equations \eqref{QNLEqulibriumEq2}
 about the uniform deformation
$\mathbf{y}_{F}.$  The linearized atomistic equation is
\begin{equation}\label{linatom}
-\la\delta^{2}\mathcal{E}^{a}\left(\mathbf{y}_{F}\right)\mathbf{u}^{a},\mathbf{w}\ra
=\la\delta \F(\mathbf{y}_F),\mathbf{w}\ra\quad\text{for all }\mathbf{w}\in \mathcal{U},
\end{equation}
and the linearized EAM-QNL equation is
\begin{equation}\label{linqnl}
-\la\delta^{2}\mathcal{E}^{qnl}\left(\mathbf{y}_{F}\right)\mathbf{u}^{qnl},\mathbf{w}\ra
=\la\delta \F(\mathbf{y}_F),\mathbf{w}\ra\quad\text{for all }\mathbf{w}\in \mathcal{U}.
\end{equation}
We thus analyze the linearized error equation
\begin{equation}\label{errorequation}
\la\delta^{2}\mathcal{E}^{qnl}\left(\mathbf{y}_{F}\right)\mathbf{e}^{qnl},\mathbf{w}\ra
=\la\mathbf{T}^{qnl},\mathbf{w}\ra\quad\text{for all }\mathbf{w}\in
\mathcal{U},
\end{equation}
where the linearized consistency error is given by
\begin{align}
\la\mathbf{T}^{qnl},\mathbf{w}\ra&:=
\la\delta^{2}\mathcal{E}^{qnl}\left(\mathbf{y}_{F}\right)
 \mathbf{u}^{a},\mathbf{w}\ra
-\la\delta^{2}\mathcal{E}^{a}\left(\mathbf{y}_{F}\right)
\mathbf{u}^{a},\mathbf{w}\ra\notag\\
& = \la\delta^{2}\hat{\mathcal{E}}^{qnl}\left(\mathbf{y}_{F}\right)
 \mathbf{u}^{a},\mathbf{w}\ra
-\la\delta^{2}\hat{\mathcal{E}}^{a}\left(\mathbf{y}_{F}\right)
\mathbf{u}^{a},\mathbf{w}\ra
\label{QNLConvergenceEq2}\\&\qquad+
\la\delta^{2}\tilde{\mathcal{E}}^{qnl}\left(\mathbf{y}_{F}\right)
 \mathbf{u}^{a},\mathbf{w}\ra
-\la\delta^{2}\tilde{\mathcal{E}}^{a}\left(\mathbf{y}_{F}\right)
\mathbf{u}^{a},\mathbf{w}\ra.
\notag
\end{align}

Now we will give an estimate of the consistency error $\mathbf{T}^{qnl},\mathbf{w}\ra$ in the following theorem.
We first define
\[\|\mathbf{v}\|^2_{\ell_\epsilon^2(\mathcal{C})}:=
\epsilon\sum_{\ell\in\mathcal{C}}v_{\ell}^2,\quad
\|\mathbf{v}\|^2_{\ell_\epsilon^2(\mathcal{I})}:=
\epsilon\sum_{\ell\in\mathcal{I}}v_{\ell}^2,\quad\text{and}\quad
\|\mathbf{v}\|^2_{\ell_\epsilon^\infty(\mathcal{I})}:=\operatorname{max}_{\ell\in\mathcal{I}}
|v_\ell|,\quad\text{for }
\mathbf{v}\in\mathcal{U},
\]
where $\mathcal{C}$ denotes the continuum region $
\{-N+1,\dots,-(K+1)\}\bigcup\{K+1,\dots,N\}$ and $\mathcal{I}$
denotes the interface $\{-(K+7),\dots, -K\}\bigcup\{K,\dots, K+7\}$.

\begin{theorem}\label{TruncationError}
The consistency error $\mathbf{T}^{qnl},\mathbf{w}\ra,$ given in \eqref{QNLConvergenceEq2}, satisfies the
following negative norm estimate
\begin{align*}
  &\left|\la\mathbf{T}^{qnl},\mathbf{w}\ra\right|
\le   \left\{ \epsilon^2\big[G''_{F}\cdot \big( (\rho'_{F})^2+12\rho'_{F}\rho'_{2F}+20(\rho'_{2F})^2\big)
-2G'_{F}\cdot\rho''_{2F}+|\phi''_{2F}|\big]\cdot \|D^{(3)}\mathbf{u}^{a}\|_{\ell^{2}_{\epsilon}(\mathcal{C})}\right.\\
&
\qquad\qquad\qquad\qquad\left.+\epsilon^{3/2}\left(C_{1}+C_{2}\right)
 \|D^{(2)}\mathbf{u}^{a}\|_{\ell^{\infty}_{\epsilon}(\mathcal{I})}
\right\}
\|D\mathbf{w}\|_{\ell^{2}_{\epsilon}}
\qquad\text{for all } \mathbf{w}\in\mathcal{U}.
\end{align*}
\end{theorem}
 \textbf{Proof}
 We focus on the first term of \eqref{QNLConvergenceEq2}
\begin{align*}
 \la\delta^{2}\hat{\mathcal{E}}^{qnl}\left(\mathbf{y}_{F}\right)
 \mathbf{u}^{a},\mathbf{w}\ra
-\la\delta^{2}\hat{\mathcal{E}}^{a}\left(\mathbf{y}_{F}\right)
\mathbf{u}^{a},\mathbf{w}\ra =
\dots+\mathbf{I}_0+\mathbf{I}_{1}+\mathbf{I}_2+\mathbf{I}_{3},
\end{align*}
 where $\mathbf{I}_0$ is associated with $\ell=0,\dots,K$,\,
$\mathbf{I}_1$ is associated with $\ell=K+1$,\, $\mathbf{I}_2$ is
associated with $\ell=K+2$ and $\mathbf{I}_3$ is associated with
$\ell=K+3,\dots,N$.

We first compute $\mathbf{I}_3$. Note that $\mathbf{u}^a$ and
$\mathbf{w}$ are $2N$-periodic, so in the calculation, when the
indices $\ell+i>N, i=1,2$, we can move these terms to the
$\{-N+1,\dots,-1\}$ part by using the periodicity as done in (6.9)
in \cite{LuskinXingjie}.  Hence, we can rearrange the terms in
$\mathbf{I}_3$ to get
\begin{align}\label{QNLConvergenceEq3}
\mathbf{I}_3 &=\epsilon\sum_{\ell=K+5}^{N}G''_{F}\cdot
\left(\rho'_{F}\right)^2
\left(-Du^{a}_{\ell-1}+2Du^{a}_{\ell}-Du^{a}_{\ell+1}\right)Dw_{\ell}\\
&\quad+\epsilon\sum_{\ell=K+5}^{N} G''_{F} \cdot
\left(\rho'_{F}\rho'_{2F}\right)
\left[4\left(-Du^{a}_{\ell-1}+2Du^{a}_{\ell}-Du^{a}_{\ell+1}\right)
+2\left(-Du^{a}_{\ell-2}+2Du^{a}_{\ell}-Du^{a}_{\ell+2}\right)\right]Dw_{\ell}\notag\\
&\quad+\epsilon\sum_{\ell=K+5}^{N} G''_{F}\cdot
\left(\rho'_{2F}\right)^2
\left[3\left(-Du^{a}_{\ell-1}+2Du^{a}_{\ell}-Du^{a}_{\ell+1}\right)
+2\left(-Du^{a}_{\ell-2}+2Du^{a}_{\ell}-Du^{a}_{\ell+2}\right)\right.\notag\\
&\quad\qquad\qquad\qquad\left.
+ \left(-Du^{a}_{\ell-3}+2Du^{a}_{\ell}-Du^{a}_{\ell+3}\right)\right]Dw_{\ell}\notag\\
&\quad+\epsilon\sum_{\ell=K+5}^{N}2G'_{F}\cdot\rho''_{2F}
\left(-Du^{a}_{\ell-1}+2Du^{a}_{\ell}-Du^{a}_{\ell+1}\right)Dw_{\ell}+
\mathbf{I}_{31}\notag
\end{align}
where $\mathbf{I}_{31}$ consists of the interfacial terms, i.e., $\ell\in\{K,\dots,K+7\},$ and
is given by the following expression
\begin{align*}
\mathbf{I}_{31}&=\epsilon G''_{F}\big\{
\left(\rho'_{F}\right)^2\left[\left(Du^a_{K+3}-Du^a_{K+4}\right)w'_{K+3}
+\left(-Du^a_{K+3}+2Du^a_{K+4}-Du^a_{K+5}\right)w'_{K+4}\right]\big.\\
&\qquad\quad+\rho'_{F}\rho'_{2F}\left[-\left(Du^a_{K+3}+Du^a_{K+4}\right)w'_{K+2}
+ \left(6Du^a_{K+3}-Du^a_{K+2}-3Du^a_{K+4}-2Du^a_{K+5}\right)w'_{K+3}\right.\big.\\
&\qquad\qquad\qquad\qquad\quad \big.\left.+ \left(12Du^a_{K+4}-Du^a_{K+2}-3Du^a_{K+3}-4Du^a_{K+5}-2Du^a_{K+6}\right)w'_{K+4} \right]\big.\\
&\quad\qquad \big. +\left(\rho'_{2F}\right)^2 \left[-\left(Du^a_{K+2}+Du^a_{K+3}+Du^a_{K+4}+Du^a_{K+5}\right)w'_{K+2}\right.\\
&\qquad\quad\qquad\left.+\left(6Du^a_{K+3}-Du^a_{K+2}-2Du^a_{K+4}-2Du^a_{K+5}-Du^a_{K+6}\right)w'_{K+3}\right.\\
&\qquad\quad\qquad \left. +\left(13Du^a_{K+4}-Du^a_{K+2}-2Du^a_{K+3}-3Du^a_{K+5}-2Du^a_{K+6}-Du^a_{K+7}\right)w'_{K+4}\right]\big\}\\
&\quad +\epsilon G'_{F}\rho''_{2F}\big\{-\left(Du^a_{K+2}+Du^a_{K+3}\right)w'_{K+2}
+\left(2Du^a_{K+3}-Du^a_{K+2}-Du^a_{K+4}\right)w'_{K+3}\\
&\qquad\qquad\qquad+\left(5Du^a_{K+4}-Du^a_{K+3}-2Du^a_{K+5}\right)w'_{K+4}\big\}.
\end{align*}
Since $\mathbf{I}_0$ is associated with $\ell=0,\dots, K$
where the QNL and the atomistic models coincide with each other, we have $\mathbf{I}_0=0$.
Similarly, by direct computation we get the following expression for the sum of $\mathbf{I}_1$ and $\mathbf{I}_2$
\begin{align}
  \mathbf{I}_1+\mathbf{I}_{2}&=
\epsilon G''_{F}\big\{ \left(\rho'_{F}\right)^2\left[\left(Du^a_{K+1}-Du^a_{K+2}\right)(w'_{K+1}-w'_{K+2})
                     \right]\big.\notag\\
&\qquad\quad\big. +\rho'_{F}\rho'_{2F}\left[\left(Du^a_{K+1}-Du^a_{K+2}\right)w'_{K}
          +\left(2Du^a_{K+1}-2Du^a_{K+2}+Du^a_{K}-Du^a_{K+3}\right)w'_{K+1}\right.\big.\notag\\
&\qquad\qquad\quad\big.\left.
             +\left(6Du^a_{K+2}-Du^a_{K}-2Du^a_{K+1}-Du^a_{K+3}\right)w'_{K+2}-\left(Du^a_{K+1}+Du^a_{K+2}\right)w'_{K+3}\right]\big.\notag\\
&\qquad \quad \big.+\left(\rho'_{2F}\right)^2\left[\left(Du^a_{K}+Du^a_{K+1}-Du^a_{K+2}-Du^a_{K+3}\right)\left(w'_{K}+w'_{K+1}\right)\right.\big.\notag\\
&\qquad\qquad \quad\big.\left.+ \left(7Du^a_{K+2}-Du^a_{K}-Du^a_{K+1}-Du^a_{K+3}\right)w'_{K+2}\right.\big.\notag\\
&\qquad\qquad\qquad\qquad \big.\left.-\left(Du^a_{K}+Du^a_{K+1}+Du^a_{K+2}+Du^a_{K+3}\right)w'_{K+3}  \right]\big\}\notag\\
&\qquad + \epsilon G'_{F}\rho''_{2F}\big\{\left(3Du^a_{K+2}-Du^a_{K+3}\right)w'_{K+2}-\left(Du^a_{K+2}+Du^a_{K+3}\right)w'_{K+3}\big\}\notag\\
&\qquad+\epsilon G''_{F}\big\{ \left(\rho'_{F}\right)^2\left[\left(Du^a_{K+2}-Du^a_{K+3}\right)(w'_{K+2}-w'_{K+3})
                     \right]\big.\label{QNLConvergenceEq4}\\
&\qquad\quad\big. +\rho'_{F}\rho'_{2F}\left[\left(Du^a_{K+2}-Du^a_{K+3}\right)w'_{K+1}
          +\left(2Du^a_{K+2}-2Du^a_{K+3}+Du^a_{K+1}-Du^a_{K+4}\right)w'_{K+2}\right.\big.\notag\\
&\qquad\qquad\quad\big.\left.
             +\left(6Du^a_{K+3}-Du^a_{K+1}-2Du^a_{K+2}-Du^a_{K+4}\right)w'_{K+3}-\left(Du^a_{K+2}+Du^a_{K+3}\right)w'_{K+4}\right]\big.\notag\\
&\qquad \quad \big.+\left(\rho'_{2F}\right)^2\left[\left(Du^a_{K+1}+Du^a_{K+2}-Du^a_{K+3}-Du^a_{K+4}\right)\left(w'_{K+1}+w'_{K+2}\right)\right.\big.\notag\\
&\qquad\qquad \quad\big.\left.+ \left(7Du^a_{K+3}-Du^a_{K+1}-Du^a_{K+2}-Du^a_{K+4}\right)w'_{K+3}\right.\big.\notag\\
&\qquad\qquad\qquad\qquad \big.\left.-\left(Du^a_{K+1}+Du^a_{K+2}+Du^a_{K+3}+Du^a_{K+4}\right)w'_{K+4}  \right]\big\}\notag\\
&\qquad + \epsilon G'_{F}\rho''_{2F}\big\{\left(3Du^a_{K+3}-Du^a_{K+4}\right)w'_{K+3}-\left(Du^a_{K+3}+Du^a_{K+4}\right)w'_{K+4}\big\}\notag.
\end{align}
Note that we can rewrite the second term of the second line of $\mathbf{I}_{3}$ as
\begin{align*}
2\left(-Du^{a}_{\ell-2}+2Du^{a}_{\ell}-Du^{a}_{\ell+2}\right)
&=2\left(-Du_{\ell-2}^{a}+2Du_{\ell-1}^a-Du_{\ell}^a\right)
+4\left(-Du_{\ell-1}^{a}+2Du_{\ell}^a-Du_{\ell+1}^a\right)\\
&\qquad\quad +2\left(-Du_{\ell}^{a}+2Du_{\ell+1}^a-Du_{\ell+2}^a\right).
\end{align*}
Similarly, we can rewrite the third term of the third line of $\mathbf{I}_{3}$ as
\begin{align*}
\left(-Du^{a}_{\ell-3}+2Du^{a}_{\ell}-Du^{a}_{\ell+3}\right)
&=\left(-Du^{a}_{\ell-3}+2Du^{a}_{\ell-2}-Du^{a}_{\ell-1}\right)
+2\left(-Du^{a}_{\ell-2}+2Du^{a}_{\ell-1}-Du^{a}_{\ell}\right)\\
&\quad +3\left(-Du^{a}_{\ell-1}+2Du^{a}_{\ell}-Du^{a}_{\ell+1}\right)
+2\left(-Du^{a}_{\ell}+2Du^{a}_{\ell+1}-Du^{a}_{\ell+2}\right)\\
&\qquad\quad +\left(-Du^{a}_{\ell+1}+2Du^{a}_{\ell+2}-Du^{a}_{\ell+3}\right).
\end{align*}
Then we combine $\mathbf{I}_{1}$, $\mathbf{I}_2$ and $\mathbf{I}_3$
together and rearrange the interfacial terms, i.e.,
$\ell\in\{K,\dots,K+7\}$. We find that the coefficients of the
interfacial terms $\mathbf{I}_{1}+\mathbf{I}_2+\mathbf{I}_{31}$ are
perfectly matched so that they are of order $\epsilon$, thus we
 obtain the following estimate by the
Cauchy-Schwarz inequality:
\begin{align}\label{QNLConvergenceEq5}
\begin{split}
&\left|
 \la\delta^{2}\hat{\mathcal{E}}^{qnl}\left(\mathbf{y}_{F}\right)
 \mathbf{u}^{a},\mathbf{w}\ra
-\la\delta^{2}\hat{\mathcal{E}}^{a}\left(\mathbf{y}_{F}\right)
\mathbf{u}^{a},\mathbf{w}\ra \right|\\
 & \quad\le \left\{\left[G''_{F}\cdot \big( (\rho'_{F})^2+12\rho'_{F}\rho'_{2F}+20(\rho'_{2F})^2\big)
-2G'_{F}\cdot\rho''_{2F}\right]\epsilon^2\cdot \|D^{(3)}\mathbf{u}^{a}\|_{\ell_\epsilon^{2}(\mathcal{C})}\right.\\
&\qquad\qquad \left.+C_{1}\epsilon\cdot \|D^{(2)}\mathbf{u}^{a}\|_{\ell_\epsilon^{2}(\mathcal{I})}\right\} \|D\mathbf{w}\|_{\ell_\epsilon^{2}}\\
\end{split}
\end{align}
where $\mathcal{I}$ is the interface: $\{K,\dots,K+7\}$,
 and $C_1$ is a constant
independent of $\epsilon.$
We note that
\begin{align*}
 \|D^{2}\mathbf{u}^{a}\|^{2}_{\ell_\epsilon^{2}(\mathcal{I})}
&=\epsilon\sum_{\ell=K}^{K+7}\left|D^{(2)}u^{a}_{\ell}\right|^{2}
\le \|D^{(2)}\mathbf{u}^{a}\|^{2}_{\ell^{\infty}_{\epsilon}(\mathcal{I})}\sum_{\ell=K}^{K+7}\epsilon
=8\epsilon
\|D^{(2)}\mathbf{u}^{a}\|^{2}_{\ell^{\infty}_{\epsilon}(\mathcal{I})}.
\end{align*}
Thus, we obtain
\begin{align*}
&\left|
 \la\delta^{2}\hat{\mathcal{E}}^{qnl}\left(\mathbf{y}_{F}\right)
 \mathbf{u}^{a},\mathbf{w}\ra
-\la\delta^{2}\hat{\mathcal{E}}^{a}\left(\mathbf{y}_{F}\right)
\mathbf{u}^{a},\mathbf{w}\ra \right|\\
&\quad \le   \left\{ \epsilon^2\big[G''_{F}\cdot \big(
(\rho'_{F})^2+12\rho'_{F}\rho'_{2F}+20(\rho'_{2F})^2\big)
-2G'_{F}\cdot\rho''_{2F}\big]\cdot
\|D^{(3)}\mathbf{u}^{a}\|_{\ell^{2}_{\epsilon}(\mathcal{C})}\right.\\
&\qquad\qquad\qquad\left.+\epsilon^{3/2}C_{1}
\|D^{(2)}\mathbf{u}^{a}\|_{\ell^{\infty}_{\epsilon}(\mathcal{I})}\right\}\cdot
\|D\mathbf{w}\|_{\ell^{2}_{\epsilon}}.
\end{align*}

We can estimate the pair potential consistency error,
$
 \la\delta^{2}\tilde{\mathcal{E}}^{qnl}\left(\mathbf{y}_{F}\right)
 \mathbf{u}^{a},\mathbf{w}\ra
-\la\delta^{2}\tilde{\mathcal{E}}^{a}\left(\mathbf{y}_{F}\right)
\mathbf{u}^{a},\mathbf{w}\ra,$ by considering the
above estimate for an
embedding energy  $G(\tilde\phi)=\tilde\phi/2$ to obtain
\begin{align*}
&\left|
 \la\delta^{2}\tilde{\mathcal{E}}^{qnl}\left(\mathbf{y}_{F}\right)
 \mathbf{u}^{a},\mathbf{w}\ra
-\la\delta^{2}\tilde{\mathcal{E}}^{a}\left(\mathbf{y}_{F}\right)
\mathbf{u}^{a},\mathbf{w}\ra \right|\\
&\quad \le \left\{ \epsilon^2 |\phi''_{2F}|
\|D^{(3)}\mathbf{u}\|_{\ell^{2}_{\epsilon}(\mathcal{C})}
+C_{2}\epsilon
\|D^{(2)}\mathbf{u}^{a}\|_{\ell^{2}_{\epsilon}(\mathcal{I})}\right\}\|D\mathbf{w}\|_{\ell^{2}_{\epsilon}} \\
&\quad \le \left\{ \epsilon^2 |\phi''_{2F}|
\|D^{(3)}\mathbf{u}\|_{\ell^{2}_{\epsilon}(\mathcal{C})}
+C_{2}\epsilon^{3/2}
\|D^{(2)}\mathbf{u}^{a}\|_{\ell^{\infty}_{\epsilon}(\mathcal{I})}\right\}\|D\mathbf{w}\|_{\ell^{2}_{\epsilon}}.
\end{align*}

 Therefore, we obtain the following optimal order estimate for the consistency error
\eqref{QNLConvergenceEq2}
\begin{align}\label{QNLConvergenceEq6}
 \begin{split}
  \left|\la\mathbf{T}^{qnl},\mathbf{w}\ra\right|
&\le \left|\la\delta^{2}\hat{\mathcal{E}}^{qnl}\left(\mathbf{y}_{F}\right)
 \mathbf{u}^{a},\mathbf{w}\ra
-\la\delta^{2}\hat{\mathcal{E}}^{a}\left(\mathbf{y}_{F}\right)
\mathbf{u}^{a},\mathbf{w}\ra \right|\\&\qquad\qquad+
\left|\la\delta^{2}\tilde{\mathcal{E}}^{qnl}\left(\mathbf{y}_{F}\right)
 \mathbf{u}^{a},\mathbf{w}\ra
-\la\delta^{2}\tilde{\mathcal{E}}^{a}\left(\mathbf{y}_{F}\right)
\mathbf{u}^{a},\mathbf{w}\ra \right|\notag\\
&\le   \left\{ \epsilon^2\big[G''_{F}\cdot \big( (\rho'_{F})^2+12\rho'_{F}\rho'_{2F}+20(\rho'_{2F})^2\big)
-2G'_{F}\cdot\rho''_{2F}+|\phi''_{2F}|\big]\cdot \|D^{(3)}\mathbf{u}^{a}\|_{\ell^{2}_{\epsilon}(\mathcal{C})}\right.\\
&
\qquad\qquad\qquad\qquad\left.+\epsilon^{3/2}\left(C_{1}+C_{2}\right)
 \|D^{(2)}\mathbf{u}^{a}\|_{\ell^{\infty}_{\epsilon}(\mathcal{I})}
\right\}
\|D\mathbf{w}\|_{\ell^{2}_{\epsilon}}
\qquad\text{for all } \mathbf{w}\in\mathcal{U}.
\end{split}
\end{align}
\quad \qed

We can now give the convergence result for the linearized EAM-QNL model.
\begin{theorem}\label{QNLconvergenceThm}
 Suppose that $\hat{A}_{F}+\tilde{A}_{F}>0$, where $\hat{A}_{F}$ and $\tilde{A}_{F}$
 are defined in \eqref{EAMCond1} and \eqref{PairAtomCond1}, and that
 \eqref{EAMfunAssumption1} and \eqref{EAMfunAssumption2} holds. Then the linearized
atomistic problem \eqref{linatom} as well as the linearized QNL
approximation \eqref{linqnl} have unique solutions, and
they satisfy the error estimate
\begin{align*}
& \|D\mathbf{y}^{a}-D\mathbf{y}^{qnl}\|_{\ell^{2}_{\epsilon}}=
\|D\mathbf{u}^{a}-D\mathbf{u}^{qnl}\|_{\ell^{2}_{\epsilon}}\\
&\,\qquad\le\frac{
 \epsilon^2\big[G''_{F}\cdot \big( (\rho'_{F})^2+12\rho'_{F}\rho'_{2F}+20(\rho'_{2F})^2\big)
-2G'_{F}\cdot\rho''_{2F}+|\phi''_{2F}|\big]\cdot \|D^{(3)}\mathbf{u}^{a}\|_{\ell^{2}_{\epsilon}(\mathcal{C})}
}
{\hat{A}_{F}+\tilde{A}_{F}}\\
&\qquad\qquad+\frac{\epsilon^{3/2}\left(C_{1}+C_{2}\right)
 \|D^{(2)}\mathbf{u}^{a}\|_{\ell^{\infty}_{\epsilon}(\mathcal{I})}}{\hat{A}_{F}+\tilde{A}_{F}}.
\end{align*}
\end{theorem}
\textbf{Proof.}
The error estimate for the EAM-QNL model follows from the error equation \eqref{errorequation},
the stability estimate in Theorem \ref{QNLStabThm},
and the consistency estimate in Theorem \ref{TruncationError}.
\quad\qed

\section{Conclusion.}
We describe a one-dimensional QNL method for the EAM potential
following~\cite{Shimokawa:2004}, and
we study the stability and convergence of a linearization of the
 next-nearest neighbor EAM-QNL energy.
 We identify conditions for the pair potential, electron density
function, and embedding function so that the lattice stability of the atomistic and the
EAM-QNL models are asymptotically equal.  These condition are necessary to guarantee
that
$u_{\ell}'=\sin(\epsilon \ell \pi)$ is the eigenfunction corresponding to
the smallest eigenvalue of $\la\delta^{2}\mathcal{E}^{a}(\mathbf{y}_{F})\mathbf{u},\mathbf{u}\ra$
with respect to the norm $\|D\mathbf{u}\|_{\ell_{\epsilon}^{2}}.$

  We then give a negative norm estimate for the consistency error and
  generalize the conclusions in \cite{Dobson:2008b} to the EAM case.
  We compare the equilibria of the atomistic and
EAM-QNL models and give an  optimal order
O($\epsilon^{3/2}$) error
  estimate for the $\ell_{\epsilon}^{2}$ norm of the
  strain in terms of the deformation in the continuum region.

\section{Acknowledgements}
We appreciate the help from Dr. Christoph Ortner and Brian Van
Koten.


\end{document}